\documentclass[12pt,reqno]{article}

\usepackage[usenames]{color}
\usepackage{amssymb}
\usepackage{amsmath}
\usepackage{amsthm}
\usepackage{amsfonts}
\usepackage{amscd}
\usepackage{graphicx}

\usepackage[colorlinks=true,
linkcolor=webgreen,
filecolor=webbrown,
citecolor=webgreen]{hyperref}

\definecolor{webgreen}{rgb}{0,.5,0}
\definecolor{webbrown}{rgb}{.6,0,0}
\usepackage{color}
\usepackage{fullpage}
\usepackage{float}

\usepackage{cases}

\usepackage{graphics,amsmath,amssymb}
\usepackage{amsthm}
\usepackage{amsfonts}
\usepackage{latexsym}

\begin{document}

\theoremstyle{plain}
\newtheorem{theorem}{Theorem}
\newtheorem{corollary}[theorem]{Corollary}
\newtheorem{lemma}[theorem]{Lemma}
\newtheorem{proposition}[theorem]{Proposition}

\theoremstyle{definition}
\newtheorem{definition}[theorem]{Definition}
\newtheorem{example}[theorem]{Example}
\newtheorem{conjecture}[theorem]{Conjecture}

\theoremstyle{remark}
\newtheorem{remark}[theorem]{Remark}

\begin{center}
\vskip 1cm{\LARGE\bf Figurate Numbers and Sums of Powers \\
\vskip .1in of Integers}
\vskip 1cm \large Jos\'{e} Luis Cereceda \\
Collado Villalba, 28400 (Madrid), Spain \\
\href{mailto:jl.cereceda@movistar.es}{\tt jl.cereceda@movistar.es}
\end{center}

\vskip .2 in

\begin{abstract}
Recently, Marko and Litvinov (ML) conjectured that, for all positive integers $n$ and $p$, the $p$-th power of $n$ admits the representation $n^p = \sum_{\ell =0}^{p-1} (-1)^{l} c_{p,\ell} F_{n}^{p-\ell}$, where $F_{n}^{p-\ell}$ is the $n$-th hyper-tetrahedron number of dimension $p-\ell$ and $c_{p,\ell}$ denotes the number of $(p -\ell)$-dimensional facets formed by cutting the $p$-dimensional cube $0 \leq x_1, x_2, \ldots, x_p \leq n-1$. In this paper we show that the ML conjecture is true for every natural number $p$. Our proof relies on the fact that the validity of the ML conjecture necessarily implies that $c_{p,\ell} = (p-\ell)! S(p, p-\ell)$, where $S(p,p-\ell)$ are the Stirling numbers of the second kind. Furthermore, we provide a number of equivalent formulas expressing the sum of powers $\sum_{i=1}^{n} i^p$ as a linear combination of figurate numbers.
\end{abstract}

\section{Introduction}
\label{sec:1}

Let $F_{n}^{k} = \binom{n+k-1}{k}$ be the $n$-th hyper-tetrahedron number of dimension $k$. In particular, $F_{n}^{2}$ is the $n$-th triangular number $T_n = \frac{1}{2}n(n+1)$ and $F_{n}^{3}$ is the $n$-th tetrahedral number $Te_n = \sum_{i=1}^{n} T_i = \frac{1}{6} n(n+1)(n+2)$ (see, e.g., \cite[Chapter 2]{guy}). Recently, Marko and Litvinov conjectured (see \cite[Conjecture 16]{marko}) that, for all positive integers $n$ and $p$, the $p$-th power of $n$ can be put as
\begin{equation}\label{ml1}
n^p = \sum_{\ell=0}^{p-1} (-1)^{\ell} c_{p,\ell} F_{n}^{p-\ell},
\end{equation}
for certain positive integer coefficients $c_{p,0}, c_{p,1},\ldots ,c_{p,p-1}$ (note the corrected factor $(-1)^\ell$ in \eqref{ml1} instead of the original one $(-1)^{p-\ell}$ appearing in \cite{marko}). Specifically, $c_{p,\ell}$ is the number of $(p-\ell)$-dimensional simplices defined by $0 \leq x_1, x_2, \ldots, x_p \leq n-1$ in conjunction with the conditions
\begin{equation*}
x_{\sigma_1} L_1 x_{\sigma_2} L_2 \ldots  L_{p-1} x_{\sigma_p},
\end{equation*}
where exactly $\ell$ symbols $L_i$ are ``$=$'', the remaining $p-\ell-1$ symbols $L_i$ are ``$\geq$'', and where $\sigma$ is a permutation of $\{ 1, 2, \ldots, p \}$. As indicated by Marko and Litvinov \cite[p.\ 18]{marko}, every such simplex is then a $(p-\ell)$-dimensional facet of the $p$-dimensional simplex $x_{\sigma_1} \geq  x_{\sigma_2} \geq \cdots \geq x_{\sigma_p}$ that is formed by cutting the original $p$-dimensional cube $0 \leq x_1, x_2, \ldots, x_p \leq n-1$.

An alternative characterization of the coefficients $c_{p,\ell}$ can be made in terms of $m$-tuples $(k_1,k_2,\ldots, k_m)$ of nonnegative integers with content $\ell = \sum_{i=1}^{m} k_i$ and support $s$ (the latter being defined as the number of indices $i$ such that $k_i >0$). As shown in \cite[Proposition 14]{marko}, the coefficients $c_{p,\ell}$ are given by the combinatorial formula
\begin{equation}\label{eq2}
c_{p,\ell} = \sum_{(k_1,k_2,\ldots,k_m)} \frac{p!}{(k_1+1)! (k_2 +1)! \cdots (k_m+1)!},
\end{equation}
where the sum runs over all $m$-tuples of nonnegative integers having the content $\ell$ and the support $s = m + \ell + 1 - p$, and such that $k_j >0$ implies that $k_{j+1} =0$ for every $j < m$.

\begin{example}
As a simple example, we may use formula \eqref{eq2} to calculate the coefficients $c_{p,\ell}$ for $p=5$. As is readily verified, for this case the list
of $m$-tuples satisfying the above conditions is: $(0,0,0,0)$ for $\ell=0$; $(1,0,0,0)$, $(0,1,0,0)$, $(0,0,1,0)$, and $(0,0,0,1)$ for $\ell=1$; $(2,0,0)$, $(0,2,0)$, $(0,0,2)$, $(1,0,1,0)$, $(0,1,0,1)$, and $(1,0,0,1)$ for $\ell=2$; $(3,0)$, $(0,3)$, $(1,0,2)$, and $(2,0,1)$ for $\ell=3$; and 4 for $\ell=4$. Hence, from \eqref{eq2} we get $c_{5,0} =120$, $c_{5,1}=240$, $c_{5,2}=150$, $c_{5,3}=30$, and $c_{5,4}=1$.
\end{example}

Furthermore, by using \eqref{eq2}, one can also deduce that, for example, for all $p \geq 1$, $c_{p,0} =p!$; for all $p \geq 2$, $c_{p,1} = \frac{1}{2} (p-1) p!$; for all $p \geq 3$, $c_{p,2} = \frac{1}{24}p!(p-2)(3p - 5)$; and that, for all $p \geq 1$, $c_{p,p-1} =1$.

Consider now the sum of powers of integers $\Sigma_{n}^{p} = \sum_{i=1}^{n} i^p$. Since $F_{n}^{k} = \sum_{i=1}^{n} F_{i}^{k-1}$,
representation \eqref{ml1} for $n^p$ immediately implies that
\begin{equation}\label{ml2}
\Sigma_{n}^{p} = \sum_{\ell=0}^{p-1} (-1)^{\ell} c_{p, \ell} F_{n}^{p+1-\ell},
\end{equation}
expressing $\Sigma_{n}^{p}$ as a linear combination of figurate numbers. In what follows, we refer to either \eqref{ml1} or \eqref{ml2} (with the coefficients $c_{p,\ell}$ given by \eqref{eq2}) as the ML conjecture. The crucial point we want to remark here is that $\Sigma_{n}^{p}$ can, in fact, be expressed in the polynomial form (see, e.g., \cite[Equation (7.5)]{gould}, \cite{pour}, and \cite[Section 4]{cere1})
\begin{equation}\label{eq4}
\Sigma_{n}^{p} = \sum_{i=1}^{p} (-1)^{p-i} i! S(p,i) \binom{n+i}{i+1},
\end{equation}
where $S(p,i)$ are the Stirling numbers of the second kind. Correspondingly, formula \eqref{eq4} can be written in terms of figurate numbers as
\begin{equation}\label{eq5}
\Sigma_{n}^{p} = \sum_{\ell =0}^{p-1} (-1)^{\ell} (p-\ell)! S(p,p-\ell) F_{n}^{p+1-\ell}.
\end{equation}
For example, setting $p=5$ in \eqref{eq5} yields
\begin{equation*}
\Sigma_{n}^{5} = 120F_{n}^{6} - 240 F_{n}^{5} + 150 F_{n}^{4} - 30 F_{n}^{3} + F_{n}^{2}.
\end{equation*}
Comparing the formulas for $\Sigma_{n}^{p}$ in \eqref{ml2} and \eqref{eq5}, and noting that the polynomials representing the figurate numbers $F_{n}^{k}$ are linearly independent, it is clear that if formula \eqref{ml2} for $\Sigma_{n}^{p}$ is true then necessarily the coefficients $c_{p,\ell}$ in \eqref{eq2} should be of the form $c_{p,\ell} = (p-\ell)! S(p, p-\ell)$ for $\ell =0,1,\ldots,p-1$. Conversely, if $c_{p,\ell} = (p-\ell)! S(p, p-\ell)$ for $\ell =0,1,\ldots,p-1$, then formula \eqref{ml2} for $\Sigma_{n}^{p}$ is true by virtue of \eqref{eq5}. This can be summarized as follows.
\begin{proposition}\label{pro2}
For $p \geq 1$ and $\ell =0,1,\ldots, p-1$, let $c_{p, \ell}$ be the coefficient defined in \eqref{eq2}. Then,
\begin{equation*}
\text{ML conjecture} \,\,  \Leftrightarrow  \,\, c_{p,\ell} = (p-\ell)! S(p, p-\ell).
\end{equation*}
\end{proposition}

The rest of the paper is organized as follows. In Section \ref{sec:2}, we give an alternative derivation of Proposition \ref{pro2} by determining the elements of the transition matrix connecting the bases $\{n,n^2,\ldots,n^p \}$ and $\{ F_{n}^{1},  F_{n}^{2}, \ldots, F_{n}^{p} \}$, and its inverse. In Section \ref{sec:3}, we prove that the ML conjecture is true. Specifically, by using the representation of the Stirling numbers of the second kind given in equation \eqref{rep}, we show that the coefficients defined in \eqref{eq2} are indeed identical to $c_{p,\ell} = (p-\ell)! S(p, p-\ell)$ for all $p \geq 1$ and $\ell =0,1, \ldots, p-1$. In Section \ref{sec:4}, we give several alternative representations for the coefficients $c_{p,\ell}$. Finally, in Section \ref{sec:5}, we provide a number of equivalent formulas expressing the sum of powers $\sum_{i=1}^{n} i^p$ as a linear combination of figurate numbers.

\section{Matrix formulation}
\label{sec:2}

For $k \geq 1$, the figurate numbers $F_{n}^{k}$ can be expanded in the basis $\{n,n^2,\ldots,n^k \}$ as
\begin{equation}\label{exp1}
F_{n}^{k} =  \binom{n+k-1}{k} = \frac{n(n+1)(n+2) \cdots (n+k-1)}{k!} = \frac{1}{k!}\sum_{r=1}^{k} s(k,r) n^r,
\end{equation}
where the numbers $s(k,r)$ are (unsigned) Stirling numbers of the first kind. For example, we have that
\begin{align*}
F_{n}^1 & = n, \\
F_{n}^2 & = \tfrac{1}{2}n + \tfrac{1}{2}n^2, \\
F_{n}^3 & = \tfrac{1}{3}n + \tfrac{1}{2}n^2 + \tfrac{1}{6} n^3,   \\
F_{n}^4 & = \tfrac{1}{4}n + \tfrac{11}{24}n^2 + \tfrac{1}{4}n^3 + \tfrac{1}{24} n^4, \\
F_{n}^5 & = \tfrac{1}{5}n + \tfrac{5}{12}n^2 + \tfrac{7}{24}n^3 + \tfrac{1}{12}n^4 + \tfrac{1}{120} n^5,
\end{align*}
or, in matrix form,
\begin{equation*}
\left( \begin{array}{c}
F_{n}^1 \\ F_{n}^2 \\ F_{n}^3 \\ F_{n}^4 \\ F_{n}^5
\end{array}\right) =
\left( \begin{array}{ccccc}
1 & 0 & 0 & 0 & 0 \\[1mm]
\tfrac{1}{2} & \tfrac{1}{2} & 0 & 0 & 0 \\[1mm]
\tfrac{1}{3} & \tfrac{1}{2} & \tfrac{1}{6} & 0 & 0 \\[1mm]
\tfrac{1}{4} & \tfrac{11}{24} & \tfrac{1}{4} & \tfrac{1}{24} & 0 \\[1mm]
\tfrac{1}{5} & \tfrac{5}{12} & \tfrac{7}{24} & \tfrac{1}{12} & \tfrac{1}{120}
\end{array}\right)
\left( \begin{array}{c}
n \\ n^2 \\ n^3 \\ n^4 \\ n^5
\end{array}\right) =
A_5 \left( \begin{array}{c}
n \\ n^2 \\ n^3 \\ n^4 \\ n^5
\end{array}\right),
\end{equation*}
where, following \cite{marko}, we call the matrices $A_p$ (for any $p \geq 1$) {\it Fermat matrices}. By inverting $A_5$ we get the transition matrix from the basis $\{F_{n}^1, F_{n}^2, F_{n}^3, F_{n}^4, F_{n}^5 \}$ to $\{n, n^2, n^3, n^4, n^5 \}$, namely
\begin{equation*}
\left( \begin{array}{c}
n \\ n^2 \\ n^3 \\ n^4 \\ n^5
\end{array}\right) =
\left( \begin{array}{ccccc}
1 & 0 & 0 & 0 & 0 \\
-1 & 2 & 0 & 0 & 0 \\
1 & -6 & 6 & 0 & 0 \\
-1 & 14 & -36 & 24 & 0 \\
1 & -30 & 150 & -240 & 120
\end{array}\right)
\left( \begin{array}{c}
F_{n}^1 \\ F_{n}^2 \\ F_{n}^3 \\ F_{n}^4 \\ F_{n}^5
\end{array}\right).
\end{equation*}
Moreover, recalling that $F_{n}^{k} = \sum_{i=1}^{n} F_{i}^{k-1}$, we can express the last matrix equation in the equivalent way
\begin{equation*}
\left( \begin{array}{c}
\Sigma_{n}^{1} \\ \Sigma_{n}^{2} \\ \Sigma_{n}^{3}  \\ \Sigma_{n}^{4} \\ \Sigma_{n}^{5}
\end{array}\right) =
\left( \begin{array}{ccccc}
1 & 0 & 0 & 0 & 0 \\
-1 & 2 & 0 & 0 & 0 \\
1 & -6 & 6 & 0 & 0 \\
-1 & 14 & -36 & 24 & 0 \\
1 & -30 & 150 & -240 & 120
\end{array}\right)
\left( \begin{array}{c}
F_{n}^2 \\ F_{n}^3 \\ F_{n}^4 \\ F_{n}^5 \\ F_{n}^6
\end{array}\right).
\end{equation*}

Let $A_p$ be a Fermat matrix, and let $a_{k,j}^{\prime}$ denote the elements of $A_{p}^{-1}$, $k,j = 1,2,\ldots, p$. Thus, for arbitrary $p \geq 1$, we have
\begin{equation*}
\left( \begin{array}{c}
\Sigma_{n}^{1} \\ \Sigma_{n}^{2} \\ \vdots \\ \Sigma_{n}^{p}
\end{array}\right) =
\left( \begin{array}{cccc}
a_{1,1}^{\prime} & a_{1,2}^{\prime} & \cdots & a_{1,p}^{\prime} \\
a_{2,1}^{\prime} & a_{2,2}^{\prime} & \cdots & a_{2,p}^{\prime}  \\
\vdots & \vdots & \ddots & \vdots  \\
a_{p,1}^{\prime} & a_{p,2}^{\prime} & \cdots & a_{p,p}^{\prime}
\end{array}\right)
\left( \begin{array}{c}
F_{n}^2 \\ F_{n}^3 \\ \vdots \\  F_{n}^{p+1}
\end{array}\right),
\end{equation*}
from which it follows that
\begin{equation}\label{exp2}
\Sigma_{n}^{p} = \sum_{i=1}^{p} a_{p,i}^{\prime} F_{n}^{i+1} = \sum_{i=0}^{p-1} a_{p,p-i}^{\prime} F_{n}^{p+1-i}.
\end{equation}
Comparing the rightmost side of \eqref{exp2} with \eqref{ml2}, and noting that the polynomials $F_{n}^{2}, F_{n}^{3}, \ldots, F_{n}^{p+1}$ are linearly independent, it is concluded that the ML conjecture is equivalent to having (cf.\ \cite{marko})
\begin{equation}\label{ml3}
a_{p,i}^{\prime} = (-1)^{p-i} c_{p,p-i}, \quad i = 1,2\ldots, p.
\end{equation}

Next, we show the following result concerning the elements $a_{k,j}^{\prime}$.

\begin{proposition}\label{pro3}
The elements $a_{k,j}^{\prime}$ of the inverse of the Fermat matrix, $A_{p}^{-1}$, are given by
\begin{equation}\label{inv}
a_{k,j}^{\prime} = (-1)^{k-j} j! S(k,j), \quad k,j =1,2,\ldots,p.
\end{equation}
\end{proposition}
\begin{proof}
Let us denote the elements of the Fermat matrix $A_{p}$ as $a_{k,j}$. Then, from \eqref{exp1} it is clear that
\begin{equation}\label{fermat}
a_{k,j} = \frac{s(k,j)}{k!}, \quad k,j =1,2,\ldots,p.
\end{equation}
Now, consider the product matrices $M_p = A_p A_{p}^{-1}$ and $N_p = A_{p}^{-1} A_{p}$, where $A_{p}^{-1}$ and $A_p$ have elements
given in \eqref{inv} and \eqref{fermat}, respectively. It is easily verified that the matrix elements $M_{k,j}$ and $N_{k,j}$ of $M_p$ and $N_p$ are
\begin{align*}
M_{k,j} & = \frac{j!}{k!} \sum_{r=1}^{p} (-1)^{r-j} s(k,r) S(r,j),  \\
N_{k,j} & = \sum_{r=1}^{p} (-1)^{k-r} S(k,r) s(r,j),
\end{align*}
for $1 \leq k,j \leq p$. Therefore, invoking the well-known orthogonality relations  (see, e.g., \cite[Thm.\ 6.24]{temple})
\vspace{-1mm}
\begin{align*}
\sum_{r=0}^{k} (-1)^{r} s(k,r) S(r,j) & = (-1)^k \delta_{k,j},  \\
\sum_{r=0}^{k} (-1)^{r} S(k,r) s(r,j) & = (-1)^k \delta_{k,j},
\end{align*}
and taking into account that $s(k,0) = S(k,0) =0$ for $k \geq 1$, and that $s(k,j) = S(k,j) = 0$ for $k < j$, we obtain $M_{k,j} = (-1)^{k-j}(j!/k!) \delta_{k,j}$ and $N_{k,j} = (-1)^{2k} \delta_{k,j}$,  and thus both $M_p$ and $N_p$ turn out to be the identity matrix $I_p$.
\end{proof}

Note that, for the case in which $k =p$, equation \eqref{inv} reads (after renaming the index $j$ as $i$) $a_{p,i}^{\prime} = (-1)^{p-i} i! S(p,i)$. Hence, from \eqref{ml3}, we conclude that the ML conjecture is true if, and only if, $c_{p,p-i} = i! S(p,i)$ or, equivalently, $c_{p,\ell} = (p-\ell)! S(p,p-\ell)$, for $\ell =0,1,\ldots, p-1$, thus recovering Proposition \ref{pro2}.

\begin{table}[b]
\centering
\begin{tabular}{|l|rrrrrrrrr|}
\hline $p \backslash  \, \ell  $ & 0 & 1 & 2 & 3 & 4 & 5 & 6 & 7 & 8   \\ \hline\hline
1 &  1 &  &  &  &  &  &  &  &   \\
2 &  2 & 1 &  &  &  &  &  &  &   \\
3 &  6 & 6 & 1 &  &  &  &  &  &   \\
4 &  24 & 36 & 14 & 1 &  &  &  &  &   \\
5 &  120 & 240 & 150 & 30 & 1 &  &  &  &   \\
6 &  720 & 1800 & 1560 & 540 & 62 & 1 &  &  &   \\
7 &  5040 & 15120 & 16800 & 8400 & 1806 & 126 & 1 &  &  \\
8 &  40320 & 141120 & 191520 & 126000 & 40824 & 5796 & 254 & 1 &  \\
9 &  362880 & 1451520 & 2328480 & 1905120 & 834120 & 186480 & 18150 & 510 & 1  \\ \hline
\end{tabular}\vspace{2mm}
\caption{Triangular array for the numbers $c_{p,\ell}$ up to $p =9$.}
\label{tb:1}
\end{table}

Table \ref{tb:1} displays the first few rows of the triangular array for the numbers $c_{p,\ell} = (p-\ell)! S(p, p-\ell)$, where $\ell =0,1,\ldots,p-1$. It is worth pointing out that, starting from the well-known triangular recurrence relation for the Stirling numbers of the second kind, namely, $S(k,j) = j S(k-1,j) + S(k-1,j-1)$ (with initial conditions $S(0,0) =1$ and $S(0,j) = S(j,0) =0$ for $j >0$), one can derive the following recurrence relation which is fulfilled by the numbers $c_{p,\ell}$:
\begin{equation}\label{recur}
c_{p,\ell} = \begin{cases}
p!,  &\text{if $\ell =0$;} \\
(p - \ell ) (c_{p-1,\ell} + c_{p-1,\ell-1}) , &\text{if $0 < \ell < p-1$;} \\
1, &\text{if $\ell = p-1$.}
\end{cases}
\end{equation}
Of course, the entries in Table \ref{tb:1} can be computed using the recursive formula \eqref{recur}. For example, $c_{9,3}$ is determined by the values of $c_{8,3}$ and $c_{8,2}$ as follows: $c_{9,3} =6(c_{8,3} + c_{8,2}) = 6(126000 + 191520) = 1905120$. Moreover, the alternating sum of the entries in the $p$-th row in Table \ref{tb:1} is given by $\sum_{\ell =0}^{p-1} (-1)^{\ell} c_{p,\ell} =1$ for all $p \geq 1$. This quickly follows from \eqref{eq5} by noting that $\Sigma_{1}^{p} = 1$ for all $p \geq 1$.

\section{Proof of the ML conjecture}\label{sec:3}

Since the recursive formula given in \eqref{recur} completely defines the numbers $c_{p,\ell} = (p-\ell)! S(p, p-\ell)$, a way to prove the ML conjecture is to show that the coefficients $c_{p,\ell}$ defined in \eqref{eq2} satisfy the said recurrence \eqref{recur}. We shall not pursue this way here. Instead, we are going to prove the ML conjecture by using the following representation for the Stirling numbers of the second kind (see, e.g., \cite[Equation (2.27)]{mezo})
\begin{equation}\label{rep}
S(n,k) = \frac{n!}{k!} \sum_{(r_1, r_2, \ldots, r_k)} \frac{1}{r_{1}! r_{2}! \cdots r_{k}!},
\end{equation}
where the summation extends over all positive integer solutions of the equation $r_1 + r_2 + \cdots + r_k = n$. Using the above representation, we can express the numbers $c_{p,\ell} = (p- \ell)! S(p, p -\ell)$ as
\begin{equation}\label{rew1}
c_{p,\ell} = \sum_{(r_1, r_2, \ldots, r_{p- \ell})} \frac{p!}{r_{1}! r_{2}! \cdots r_{p -\ell}!},
\end{equation}
where each of the $(p-\ell)$-tuples $(r_1, r_2, \ldots, r_{p- \ell})$ in the summation has the content $\sum_{i=1}^{p-\ell} r_i = p$, with $r_i \geq 1$.

\begin{example}
For $\ell = p-2$ and $p \geq 2$, from \eqref{rew1} we have
\begin{equation*}
c_{p,p-2} = \sum_{r_1 + r_2 =p} \frac{p!}{r_{1}! r_{2}!} = \sum_{j=1}^{p-1} \frac{p!}{j! (p-j)!}
= \sum_{j=1}^{p-1} \binom{p}{j} = 2^p -2.
\end{equation*}
Similarly, for $\ell =p-3$ and $p \geq 3$, from \eqref{rew1} we have
\begin{equation}\label{trin}
c_{p,p-3} = \sum_{r_1 + r_2 + r_3 =p} \frac{p!}{r_{1}! r_{2}! r_{3}!} = \sum_{i+j+k=p} \binom{p}{i,j,k} - 3 - 3 \sum_{t=1}^{p-1}
\binom{p}{t} = 3^p - 3 \cdot 2^p +3 ,
\end{equation}
where $\binom{p}{i,j,k}$ are the trinomial coefficients, with $i$, $j$, and $k$ being nonnegative integers.
\end{example}

On the other hand, by renaming $k_i + 1$ as $j_i$ for $i =1,2,\ldots, m$, equation \eqref{eq2} can be rewritten as
\begin{equation}\label{rew2}
c_{p,\ell} = \sum_{(j_1,j_2,\ldots, j_m)} \frac{p!}{j_{1}! j_{2}! \cdots j_{m}!},
\end{equation}
where now the summation extends over all $m$-tuples of positive integers $(j_1,j_2,\ldots, j_m)$ with the content $\sum_{i=1}^{m} j_i = \ell + m$, where $m = p+t-\ell -1$ (with $t$ being the number of indices $i$ for which $j_i \geq 2$), and such that $j_i \geq 2$ implies that $j_{i+1}=1$ for every $i < m$.

\begin{example}
For $\ell = p-2$ and $p \geq2$, it can be seen that the allowed $m$-tuples for this case are: $(p-1,1)$, $(1,p-1)$, and $(a,1,b)$, with $a+b =p$ and $a,b \geq2$. From \eqref{rew2} we then have
\begin{equation*}
c_{p,p-2} =2p + \sum_{a+b=p} \frac{p!}{a! b!} = 2p + \sum_{j=2}^{p-2} \frac{p!}{j! (p -j)!} = 2p +
\sum_{j=2}^{p-2} \binom{p}{j} = 2^p -2.
\end{equation*}
Moreover, for $\ell = p-3$ and $p \geq 3$, the allowed $m$-tuples for this case are: $(p-2,1,1)$, $(1,p-2,1)$, $(1,1,p-2)$, $(a,1,b,1)$, $(a,1,1,b)$, $(1,a,1,b)$, and $(c,1,d,1,e)$, where $a+b = p-1$ with $a,b \geq 2$, and $c+d+e = p$ with $c,d,e \geq 2$. Thus, from \eqref{rew2} we obtain
\begin{equation*}
c_{p,p-3} = 3p(p-1) + 3 \sum_{a+b =p-1} \frac{p!}{a! b!} + \sum_{c+d+e =p} \frac{p!}{c! d! e!}.
\end{equation*}
Now, we have
\begin{equation*}
\sum_{c+d+e =p} \frac{p!}{c! d! e!} = \sum_{r_1 + r_2 + r_3 =p} \frac{p!}{r_{1}! r_{2}! r_{3}!} - 3p(p-1) - 3 \sum_{a+b =p-1}
\frac{p!}{a! b!},
\end{equation*}
where $r_1, r_2, r_3 \geq 1$. Therefore, it follows that
\begin{equation*}
c_{p,p-3} =  \sum_{r_1 + r_2 + r_3 =p} \frac{p!}{r_{1}! r_{2}! r_{3}!} = 3^p - 3 \cdot 2^p +3 ,
\end{equation*}
where we have used the previous result in \eqref{trin}. In the same manner, for $\ell = p-4$ and $p \geq 4$, starting from either \eqref{rew1} or \eqref{rew2} we find that
\begin{equation}\label{j4}
c_{p,p-4} = 4^p - 4 \cdot 3^p + 6 \cdot 2^p -4 .
\end{equation}
Indeed, since $j! S(p,j) = \sum_{r=0}^{j} (-1)^r \binom{j}{r} (j -r)^p$, the validity of the ML conjecture will retrospectively enable us to deduce that both \eqref{rew1} and \eqref{rew2} lead to the general result
\begin{equation*}
c_{p,p-j} = \sum_{r=0}^{j} (-1)^r \binom{j}{r} (j-r)^p , \quad j =1,2,\ldots,p.
\end{equation*}
\end{example}

In view of Proposition \ref{pro2}, and making use of \eqref{rew1} and \eqref{rew2}, it follows that the validity of the ML conjecture is then equivalent to the following number theoretic identity.

\begin{proposition}\label{pro6}
For all $p \geq 1$ and $\ell = 0,1,\ldots, p-1$, we have
\begin{equation}\label{equal}
\sum_{(j_1,j_2,\ldots, j_m)} \frac{p!}{j_{1}! j_{2}! \cdots j_{m}!} \,\, = \sum_{(r_1, r_2, \ldots, r_{p- \ell})} \frac{p!}{r_{1}! r_{2}! \cdots r_{p -\ell}!}.
\end{equation}
\end{proposition}

Next, we sketch the proof of Proposition \ref{pro6} and, as a consequence, of the ML conjecture.
\begin{proof}
The first step towards proving \eqref{equal} is to expand the summation in the left-hand side of \eqref{equal} over all the $m$-tuples $(j_1,j_2, \ldots, j_m)$ of positive integers fulfilling the conditions stipulated after equation \eqref{rew2}. Recalling that $t$ is the number of $j_i$'s for which $j_i \geq 2$, it is not hard to show that, for $\ell = p -j$, equation \eqref{rew2} can be decomposed as
\begin{numcases}{c_{p,p-j}=}
\displaystyle{\sum_{t=1}^{j} \binom{j}{t} \sum_{s_1 + s_2 + \cdots + s_t = p +t- j} \frac{p!}{s_{1}! s_{2}! \cdots s_{t}!}},
&\text{for $j =1,2,\ldots, p-1$;}  \label{proof1a} \\
p!,  &\text{for $j = p$,} \label{proof1b}
\end{numcases}
where $s_1, s_2, \ldots, s_t$ are all integers $\geq 2$. For example, for $j=4$, we have
\begin{align*}
c_{p,p-4} & = 4 \sum_{s_1 = p-3} \frac{p!}{s_{1}!} + 6  \sum_{s_1 + s_2 = p-2} \frac{p!}{s_{1}! s_{2}!}
 + 4 \sum_{s_1 + s_2 + s_3= p-1} \frac{p!}{s_{1}! s_{2}! s_{3}!} \\
 &  \quad \,\, + \sum_{s_1 + s_2 + s_3 +s_4= p} \frac{p!}{s_{1}! s_{2}! s_{3}! s_{4}!}.
\end{align*}
As an even more concrete example, we may use the last equation to calculate the coefficient $c_{9,5}$ as follows
\begin{align*}
c_{9,5} & = 4 \sum_{s_1 =6} \frac{9!}{s_{1}!} + 6 \sum_{s_1 + s_2 =7} \frac{9!}{s_{1}! s_{2}!} + 4 \sum_{s_1 +s_2 +s_3 =8}
\frac{9!}{s_{1}! s_{2}! s_{3}!} + \sum_{s_1 +s_2 + s_3 + s_4 =9} \frac{9!}{s_{1}! s_{2}! s_{3}! s_{4}!}  \\
& = 9! \left[ 4 \, \frac{1}{6!} + 6 \left( \frac{2}{2! 5!} + \frac{2}{3! 4!} \right) + 4 \left( \frac{3}{2! 2! 4!} + \frac{3}{2! 3! 3!} \right)
+ \frac{4}{2! 2! 2! 3!} \right] = 186480.
\end{align*}
Of course, equation \eqref{j4} reproduces this result as $186480 = 4^9 -4\cdot 3^9 + 6\cdot 2^9 -4$. Furthermore, as an aside, it is to be noted that the number of individual summands involved in \eqref{proof1a} amounts to
\begin{equation*}
N(p,j) =\sum_{t=1}^{j} \binom{j}{t} \binom{p-j-1}{t-1} = \binom{p-1}{j-1}.
\end{equation*}
For the example above for which $p =9$ and $j=4$, we have $N(9,4) = 4 \cdot 1 + 6 \cdot 4 + 4 \cdot 6 + 1 \cdot 4 =56$.

Returning to the main argument, let us, for convenience, rewrite \eqref{proof1a} in the form
\begin{equation}\label{proof2}
c_{p,p-j} = \sum_{s_1 + s_2 + \cdots + s_j = p} \frac{p!}{s_{1}! s_{2}! \cdots s_{j}!} + W(p,j),
\end{equation}
where
\begin{equation*}
W(p,j) = \sum_{t=1}^{j-1} \binom{j}{t} \sum_{s_1 + s_2 + \cdots + s_t = p +t- j} \frac{p!}{s_{1}! s_{2}! \cdots s_{t}!}.
\end{equation*}
The second (and last) crucial point is to realize that $W(p,j)$ is precisely the summation
\begin{equation*}
W(p,j) = \sum_{w_1 + w_2 + \cdots + w_j = p} \frac{p!}{w_{1}! w_{2}! \cdots w_{j}!},
\end{equation*}
where $w_1, w_2, \ldots, w_j$ are all integers $\geq 1$ and {\it at least} one of them is equal to $1$. (Note that the case $w_1 = w_2 = \ldots = w_j =1$ is excluded since, as far as equation \eqref{proof1a} is concerned, $j \leq p-1$.) For example, we have that
\begin{equation*}
W(p,4) = 4 \sum_{s_1 +1+1+1= p} \frac{p!}{s_{1}!1!1!1!} + 6  \sum_{s_1 + s_2 +1+1= p} \frac{p!}{s_{1}! s_{2}!1!1!}
 + 4 \sum_{s_1 + s_2 + s_3+1= p} \frac{p!}{s_{1}! s_{2}! s_{3}!1!}.
\end{equation*}\pagebreak

\noindent Now, for $\ell = p-j$, \eqref{rew1} reads as
\begin{equation}\label{proof3}
c_{p,p-j} = \sum_{r_1 + r_2 + \cdots + r_j = p} \frac{p!}{r_{1}! r_{2}! \cdots r_{j}!},
\end{equation}
where $r_1, r_2, \ldots, r_j$ are all integers $\geq 1$. Hence, it is clear that, by definition,
\begin{equation}\label{proof4}
 \sum_{r_1 + r_2 + \cdots + r_j = p} \frac{p!}{r_{1}! r_{2}! \cdots r_{j}!} = \sum_{s_1 + s_2 + \cdots + s_j = p} \frac{p!}{s_{1}! s_{2}! \cdots s_{j}!} + W(p,j),
\end{equation}
for $j =1,2,\ldots, p-1$. Therefore, noting that both \eqref{rew1} and \eqref{rew2} give $c_{p,0} = p!$, and taking into account \eqref{proof2}, \eqref{proof3}, and \eqref{proof4}, it follows that the identity \eqref{equal} holds for all $p \geq 1$ and $\ell = 0,1,\ldots, p-1$, and, consequently, it is concluded that the ML conjecture is true.
\end{proof}

We end this section with the following observation. It is a well-known, classical result that the number of surjections $T(m,n)$ from the set $A$ (having $m$ elements) onto the set $B$ (having $n$ elements, with $m \geq n$) is given by $n! S(m,n)$; see, e.g., \cite{hilton}. Since $c_{p,p-j} = j! S(p,j)$, this means, in particular, that the formula for $c_{p,p-j}$ in \eqref{proof1a}-\eqref{proof1b} actually gives the number of surjections $T(p,j)$ from a $p$-element set to a $j$-element set for each $j =1,2,\ldots,p$. Moreover, from the recurrence relation $T(p,j) = j (T(p-1,j) + T(p-1,j-1))$, we immediately get the recursive formula
\begin{equation*}
c_{p, p-j} = \begin{cases}
1,  &\text{if $j =1$;} \\
j ( c_{p-1,p-j} + c_{p-1,p-j-1} ) ,  &\text{if $0 < j < p-1$;} \\
p!, &\text{if $j = p$,}
\end{cases}
\end{equation*}
which is, of course, the same as that in \eqref{recur}.

\section{Alternative representations for the coefficients $c_{p,\ell}$}
\label{sec:4}

Next, we state the following alternative formula for the coefficients $c_{p,\ell}$:
\begin{equation}\label{form1}
c_{p, \ell} = \begin{cases}
p!,  &\text{for $\ell =0$;} \\
p! (p- \ell) \displaystyle{\sum_{k=1}^{\ell} \binom{p-\ell -1}{k-1}\binom{2\ell -p}{\ell-k} \frac{(k-1)!}{(k+\ell)!}} S(k+\ell, k), &\text{for $1 \leq \ell \leq p-1$.}
\end{cases}
\end{equation}
Recalling that $c_{p,\ell} = (p -\ell)! S(p, p-\ell)$, it is easily seen that the representation \eqref{form1} is a direct consequence of the identity (see
\cite[Thm.\ 2.3]{sun})
\begin{equation*}
S(p+\ell, p) = \sum_{k=0}^{\ell} \binom{\ell -p}{\ell -k} \binom{\ell +p}{\ell +k} S(k +\ell, k),
\end{equation*}
which is in turn implied by the results found in Gould's paper \cite{gould2}. We can write down the expressions for $c_{p,0}, c_{p,1}, \ldots, c_{p,7}$ obtained from \eqref{form1} as follows
\begin{equation}\label{set1}
\begin{split}
c_{p,0} & = p!,  \\
c_{p,1} & = \frac{p!}{2!} (p-1),  \\
c_{p,2} & = \frac{p!}{3!} (p-2) \left( \frac{3p-5}{4} \right),  \\
c_{p,3} & = \frac{p!}{4!} (p-3)^2 \left(\frac{p-2}{2} \right),  \\
c_{p,4} & = \frac{p!}{5!} (p-4) \left( \frac{15p^3 -150p^2 +485p -502}{48} \right),  \\
c_{p,5} & = \frac{p!}{6!} (p-5)^2 (p-4) \left(\frac{3p^2 -23p +38}{16} \right),  \\
c_{p,6} & = \frac{p!}{7!} (p-6) \left( \frac{63 p^5 - 1575p^4 + 15435p^3 - 73801p^2 + 171150p -152696}{576} \right),  \\
c_{p,7} & = \frac{p!}{8!} (p-7)^2 (p-6) \left(\frac{9p^4 - 198p^3 + 1563p^2 - 5182p +6008}{144} \right).
\end{split}
\end{equation}
In view of these and successive coefficients, it seems safe to conjecture that $(p-\ell)^2 (p-\ell+1)$ constitutes a factor of $c_{p,\ell}$ for all odd $\ell$ with $\ell \geq 3$. A table of coefficients related to those in \eqref{set1} is given in Equations (1.18)-(1.22) of Gould's book \cite{gould3}. It is to be noted, on the other hand, that \eqref{form1} may also be obtained by considering the generalized Bernoulli numbers or N\"{o}rlund polynomials $B_{k}^{(\alpha)}$ defined by the generating function
\begin{equation*}
\left( \frac{x}{e^x -1} \right)^{\alpha} = \sum_{k=0}^{\infty} B_{k}^{(\alpha)} \, \frac{x^k}{k!}, \quad |x| < 2\pi,
\end{equation*}
where the parameter $\alpha$ stands for any arbitrary (real or complex) number. It turns out that $B_{k}^{(\alpha)}$ is a polynomial in $\alpha$ of degree $k$. In particular, $B_{k}^{(1)} = B_k$ are the ordinary Bernoulli numbers. The Stirling numbers of the second kind are related to the generalized Bernoulli numbers by means of (see \cite[Equation (14.12)]{gould4})
\begin{equation*}
S(n+k,n) = \binom{n+k}{k} B_{k}^{(-n)}.
\end{equation*}
Therefore, we have that
\begin{equation}\label{form2}
c_{p, \ell} = (p-\ell)! S(p, p-\ell) = \frac{p!}{\ell !} B_{\ell}^{(-p + \ell)}.
\end{equation}
Moreover, $B_{\ell}^{(\alpha)}$ is given explicitly by Equation (15) of Srivastava and Todorov's paper \cite{todorov}
\begin{equation}\label{form3}
B_{\ell}^{(\alpha)} = \sum_{k=0}^{\ell} (-1)^k \binom{\alpha +k-1}{k} \binom{\alpha + \ell}{\ell -k}
\binom{k +\ell}{k}^{-1} S(k+\ell, k).
\end{equation}
Hence, taking equation \eqref{form3} with $\alpha \to -p +\ell$, and substituting the resulting $B_{\ell}^{(-p+\ell)}$ into \eqref{form2}, we find that
\begin{equation}\label{form4}
c_{p,\ell} = p! \sum_{k=0}^{\ell} (-1)^k \binom{-p+\ell +k-1}{k} \binom{2\ell -p}{\ell -k} \frac{k!}{(k+\ell)!} S(k+\ell, k).
\end{equation}
Finally, noting that $\binom{p -\ell}{k} = (-1)^k \binom{-p+\ell +k -1}{k}$, we conclude that the representation in \eqref{form4} is indeed the same as that in \eqref{form1}.

Furthermore, from the relation $S(p, p-\ell) = \binom{p}{\ell} B_{\ell}^{(-p +\ell)}$, it is readily seen that the Stirling numbers of the second kind of the type $S(p,p-\ell)$ are polynomials in $p$ of degree $2\ell$. These can be expressed, in particular, in the form
\begin{equation}\label{form5}
S(p, p-\ell) = \begin{cases}
1,  &\text{for $\ell =0$;} \\
\displaystyle{\sum_{j=0}^{\ell -1} a_{\ell,j} \binom{p}{2\ell - j}},  &\text{for $\ell =1,2,\ldots,p-1$,}
\end{cases}
\end{equation}
where the coefficients $a_{\ell,j}$, which were introduced by Jordan and Ward \cite{carlitz} (see also \cite{roman} and references therein), satisfy the recurrence
\begin{equation*}
a_{\ell,j} = (2\ell -j-1) a_{\ell-1,j} + (\ell -j) a_{\ell -1,j-1},
\end{equation*}
and are given explicitly by Equation (1.13) of Carlitz's paper \cite{carlitz}
\begin{equation}\label{form6}
a_{\ell,j} = \frac{1}{(\ell -j)!}\sum_{k=0}^{\ell -j} (-1)^k \binom{\ell -j}{k} \sum_{r=0}^{k} r! \binom{k}{r} \binom{2\ell -j}{r}
(\ell -j -k)^{2\ell -j -r}.
\end{equation}
From \eqref{form5} we can therefore express the coefficients $c_{p,\ell}$ as
\begin{equation}\label{form7}
c_{p, \ell} = \begin{cases}
p!,  &\text{for $\ell =0$;} \\
(p- \ell)! \displaystyle{\sum_{j=0}^{\ell -1} a_{\ell,j} \binom{p}{2\ell - j}},  &\text{for $\ell =1,2,\ldots,p-1$,}
\end{cases}
\end{equation}
with the $a_{\ell,j}$'s being given by \eqref{form6}. Table \ref{tb:2} shows the first few rows of the triangular array for the Jordan coefficients. Let us note, in particular, that $a_{\ell,0} = (2\ell -1)!!$.

\begin{table}[ttt]
\centering
\begin{tabular}{|l|rrrrrrrr|}
\hline $\ell \backslash  \, j  $ & 0 & 1 & 2 & 3 & 4 & 5 & 6 & 7   \\ \hline\hline
1 &  1 &  &  &  &  &  &  &   \\
2 &  3 & 1 &  &  &  &  &  &   \\
3 &  15 & 10 & 1 &  &  &  &  &     \\
4 &  105 & 105 & 25 & 1 &  &  &  &    \\
5 &  945 & 1260 & 490 & 56 & 1 &  &  &   \\
6 &  10395 & 17325 & 9450 & 1918 & 119 & 1 &  &   \\
7 &  135135 & 270270 & 190575 & 56980 & 6825 & 246 & 1 &  \\
8 &  2027025 & 4729725 & 4099095 & 1636635 & 302995 & 22935 & 501 & 1 \\  \hline
\end{tabular}\vspace{2mm}
\caption{Triangular array for the Jordan coefficients $a_{\ell,j}$ up to $\ell =8$.}
\label{tb:2}
\end{table}

On the other hand, for fixed $\ell$, the exponential generating function of the Stirling numbers of the second kind $S(p,\ell)$ is given by (see, e.g., \cite[Equation (2.19)]{mezo})
\begin{equation*}
\frac{(e^x -1 )^{\ell}}{\ell !} = \sum_{p=0}^{\infty} S(p,\ell) \frac{x^p}{p!},
\end{equation*}
from which one can deduce that
\begin{equation}\label{deriv}
c_{p, \ell} = (p-\ell)! S(p, p-\ell) = f_{p,\ell}^{(p)}(0),
\end{equation}
where $f_{p,\ell}^{(p)}(0) = d^p f_{p,\ell}(x)/ dx^p |_{x=0}$ is the $p$-th derivative of $f_{p,\ell}(x) = (e^x -1)^{p-\ell}$ with respect to $x$ evaluated in the limit when $x \to 0$. As a simple example, for $\ell = p-4$, we have
\begin{equation*}
f_{p,p-4}(x) = (e^x -1)^4 = e^{4x} - 4 e^{3x} +6 e^{2x} -4 e^x +1.
\end{equation*}
It is easily seen that
\begin{equation*}
\frac{d^p f_{p,p-4}(x)}{dx^p} = 4^p e^{4x} - 4 \cdot 3^p e^{3x} + 6 \cdot 2^p e^{2x} - 4 e^x,
\end{equation*}
and then, from \eqref{deriv}, it follows that $c_{p,p-4} = 4^p - 4 \cdot 3^p + 6 \cdot 2^p -4$, thus recovering \eqref{j4}.

We end this section by quoting two other possible representations for $c_{p,\ell}$, namely
\begin{equation}\label{form8}
c_{p,\ell} = (p- \ell)! \sum_{k=0}^{\ell} \binom{\ell -p}{\ell +k} \binom{p+ \ell}{\ell -k} s(k+\ell,k),
\end{equation}
and
\begin{equation}\label{form9}
c_{p, \ell} = \begin{cases}
p!,  &\text{for $\ell =0$;} \\
(p- \ell)! \displaystyle{\sum_{j=0}^{\ell -1}  \bigg< \!\!\! \genfrac{<}{>}{0pt}{}{\ell}{j} \!\!\! \bigg> \binom{p+\ell-1-j}{2\ell}},
&\text{for $\ell =1,2,\ldots,p-1$,}
\end{cases}
\end{equation}
where $s(k+\ell,k)$ are the (unsigned) Stirling numbers of the first kind, and $\big< \!\! \genfrac{<}{>}{0pt}{}{\ell}{j} \!\! \big>$ are the second-order Eulerian numbers indexed so that $0 \leq j \leq \ell-1$ (see \cite[Table 270 and Equation (6.43)]{graham}).

\section{Power sums as a linear combination of figurate numbers}
\label{sec:5}

Starting with \eqref{ml2} and using the representations for $c_{p,\ell}$ given in \eqref{form4}, \eqref{form7}, \eqref{form8}, and \eqref{form9}, we successively obtain the following formulas expressing $\Sigma_{n}^{p}$ as a linear combination of the figurate numbers $F_{n}^{2}, F_{n}^{3},
\ldots, F_{n}^{p+1}$:
\begin{equation}\label{set2}
\begin{split}
\Sigma_{n}^{p} & = p! \sum_{\ell=0}^{p-1} (-1)^{\ell} F_{n}^{p+1-\ell} \sum_{k=0}^{\ell}
\binom{p -\ell}{k} \binom{2\ell -p}{\ell -k}
\frac{k!}{(k+\ell)!} S(k+\ell, k), \\
\Sigma_{n}^{p} & = p!F_{n}^{p+1} + \sum_{j=1}^{p-1} (-1)^{j} (p-j)! F_{n}^{p+1-j} \sum_{k=0}^{j-1} a_{j,k} \binom{p}{2j-k}, \\
\Sigma_{n}^{p} & = \sum_{\ell=0}^{p-1} (-1)^{\ell} (p-\ell)! F_{n}^{p+1-\ell} \sum_{k=0}^{\ell} \binom{\ell-p}{\ell+k}\binom{p+\ell}
{\ell-k} s(k+\ell,k),  \\
\Sigma_{n}^{p} & = p!F_{n}^{p+1} + \sum_{j=1}^{p-1} (-1)^j (p-j)! F_{n}^{p+1-j} \sum_{k=0}^{j -1}  \bigg< \!\!\!
\genfrac{<}{>}{0pt}{}{j}{k} \!\!\! \bigg> \binom{p+j-1-k}{2j}.
\end{split}
\end{equation}
All four linear combinations in \eqref{set2} can be rewritten in the form $\sum_{\ell=0}^{p-1} f_{p,\ell} F_{n}^{p+1-\ell}$, where the coefficient $f_{p,\ell}$ is equal to $(-1)^{\ell}(p-\ell)! S(p,p-\ell)$ for each $\ell =0,1,\ldots,p-1$, in accordance with \eqref{eq5}.

In addition to the formulas in \eqref{eq5} or \eqref{set2}, there are several other alternative formulas expressing $\Sigma_{n}^{p}$ in terms of the figurate numbers $F_{n}^{k} = \binom{n+k-1}{k}$. For example, the following two well-known polynomial formulas for $\Sigma_{n}^{p}$ (see, e.g., \cite{cere1,shirali,tsao}):
\begin{align}
\Sigma_{n}^{p} & = \sum_{j=1}^{p} j! S(p,j) \binom{n+1}{j+1}, \label{stir} \\
\intertext{and}
\Sigma_{n}^{p} & = \sum_{j=1}^{p} \genfrac{<} {>}{0pt}{}{p}{j} \binom{n+j}{p+1}, \label{euler}
\end{align}
can equivalently be written in terms of $F_{n}^{k}$ as
\begin{align}
\Sigma_{n}^{p} & = \sum_{j=1}^{p} j! S(p,j) F_{n-j+1}^{j+1}, \label{alt1} \\
\intertext{and}
\Sigma_{n}^{p} & = \sum_{j=1}^{p} \genfrac{<} {>}{0pt}{}{p}{j} F_{n+j-p}^{p+1},  \label{alt2}
\end{align}
respectively, where $\genfrac{<} {>}{0pt}{}{p}{j}$ are the ordinary Eulerian numbers, with the initial value $\genfrac{<} {>}{0pt}{}{p}{1} =1$ for all $p \geq 1$. Note that \eqref{alt2} only involves figurate numbers of dimension $p+1$. Along with the above two formulas in \eqref{stir} and \eqref{euler}, we may quote another, not so well-known formula for $\Sigma_{n}^{p}$ which is a variant of that in \eqref{stir}, namely (see, e.g., \cite[Equation (9)]{roman} and \cite{cere2})
\begin{equation}\label{alt3}
\Sigma_{n}^{p} = \sum_{j=1}^{p+1} (j-1)! S(p+1,j) \binom{n}{j} = \sum_{j=1}^{p+1} (j-1)! S(p+1,j) F_{n-j+1}^{j}.
\end{equation}

As an example, for $p=8$, from \eqref{eq5}, \eqref{alt1}, \eqref{alt2}, and \eqref{alt3}, we obtain the equivalent polynomial representations
\begin{align*}
\Sigma_{n}^{8} & = 40320 F_{n}^{9} - 141120 F_{n}^{8} + 191520 F_{n}^{7} \\
& \quad - 126000 F_{n}^{6} + 40824 F_{n}^{5} - 5796 F_{n}^{4} + 254 F_{n}^{3} - F_{n}^{2}  \\
& = 40320 F_{n-7}^{9} + 141120 F_{n-6}^{8} + 191520 F_{n-5}^{7} \\
& \quad + 126000 F_{n-4}^{6} + 40824 F_{n-3}^{5} + 5796 F_{n-2}^{4} + 254 F_{n-1}^{3} + F_{n}^{2}  \\
& = F_{n}^{9} + 247 F_{n-1}^{9} + 4293 F_{n-2}^{9} + 15619 F_{n-3}^{9} \\
& \quad + 15619 F_{n-4}^{9} + 4293 F_{n-5}^{9} + 247 F_{n-6}^{9} + F_{n-7}^{9}  \\
& = 40320 F_{n-8}^{9} + 181440 F_{n-7}^{8} + 332640 F_{n-6}^{7} + 317520 F_{n-5}^{6} \\
& \quad + 166824 F_{n-4}^{5} + 46620 F_{n-3}^{4} + 6050 F_{n-2}^{3} + 255 F_{n-1}^{2} + F_{n}^{1}.
\end{align*}

For completeness, let us finally mention that, as is well known, the power sums $\Sigma_{n}^{p}$ can be expressed as polynomials in the triangular numbers $T_n$ (the so-called Faulhaber polynomials \cite{edwards}) as follows
\begin{align*}
\Sigma_{n}^{2k} & = \Sigma_{n}^{2} \left[ b_{k,0} + b_{k,1} T_n + b_{k,2} \big(  T_n \big)^2 + \cdots +
b_{k,k-1} \big(  T_n \big)^{k-1} \right], \\
\Sigma_{n}^{2k+1} & =  \big( T_n \big)^2 \left[ c_{k,0} + c_{k,1} T_n + c_{k,2} \big( T_n \big)^2 + \cdots +
c_{k,k-1} \big(  T_n \big)^{k-1} \right],
\end{align*}
where $b_{k,j}$ and $c_{k,j}$ are non-zero rational coefficients for $j =0,1,\ldots,k-1$ and $k \geq 1$. In particular, $\Sigma_{n}^{3} =
\big(T_n \big)^2$.

\section*{Note added}

After the completion of this work, a much shorter and direct proof of the ML conjecture was devised by Professor Franti\v{s}ek Marko himself \cite{marko2}. I believe, however, that the proof of the ML conjecture presented here still retains its interest in its own right.

\end{document}